\documentclass[11pt]{article}
      \usepackage{latexsym}
      \newtheorem{theorem}{Theorem}[section]

      \newtheorem{lemma}[theorem]{Lemma}

      \def\nn{\nonumber}
      \def\rf#1{\mbox{$(\ref{#1})$}}

      \def\be{\begin{equation}} %\be=\begin{equation}
      \def\ee{\end{equation}} %\ee=\end{equation}
      \def\beqn{\begin{eqnarray}} %\beqn=\begin{eqnarray}
      \def\eeqn{\end{eqnarray}} %\eeqn=\end{eqnarray}
      \def\beq{\begin{eqnarray*}} %\beq=\begin{eqnarray*}
      \def\eeq{\end{eqnarray*}}
      \def\proof{{\bf Proof:}\ }
      \def\mb{\mbox} %\mb=\mbox
      \def\ga{\gamma} %\ga=\gamma
      \def\la{\lambda} %\la=\lambda
      \def\ra{\rightarrow} %\ra=\rightarrow
      \def\esssup{\mathop {\rm ess\,sup}}
      
\begin{document}
      \title{\bf Large deviations for Dirichlet processes and Poisson-Dirichlet distributions with two parameters\thanks{Research supported by
      the Natural Science and Engineering Research Council of Canada}}
      \author{Shui Feng\\McMaster
      University}
      \date{\today}
      \maketitle
      \begin{abstract}
      Large deviation principles are established
      for the two-parameter Poisson-Dirichlet distribution and two-parameter Dirichlet process when parameter $\theta$ approaches infinity.
      The motivation for these results is to understand the differences in terms of large deviations between the two-parameter models and their
      one-parameter counterparts.
      New insight is obtained about the role of the second parameter $\alpha$
      through a comparison with the corresponding results for the
      one-parameter Poisson-Dirichlet distribution and Dirichlet process.

      \vspace*{.125in} \noindent {\it Keywords:} GEM representation, Poisson-Dirichlet
      distribution, two parameter Poisson-Dirichlet distribution, Dirichlet processes, large deviation principles.
      \vspace*{.125in}

      \end{abstract}
      \section{Introduction}
      \setcounter{equation}{0}

      For $\theta >0$ and $\alpha$ in $(0,1)$, let $U_k,
      k=1,2,...$, be a sequence of independent
      random variables such that $U_k$ has $Beta(1-\alpha,\theta+ k\alpha)$ distribution.
      Set
      \be \label{GEM1}
      X^{\theta,\alpha}_1 = U_1,\  X^{\theta,\alpha}_n = (1-U_1)\cdots (1-U_{n-1})U_n,\  n \geq 2.
      \ee

      Then with probability one $$\sum_{k=1}^{\infty}X^{\theta,\alpha}_k =1,$$ and the law of $(X^{\theta,\alpha}_1,X^{\theta,\alpha}_2,...)$
      is called the two-parameter GEM distribution denoted by $GEM(\theta,\alpha).$

      Let ${\bf P}(\theta,\alpha)=(P_1(\theta,\alpha), P_2(\theta,\alpha),...)$ denote $(X^{\theta,\alpha}_1,X^{\theta,\alpha}_2,...)$
      in descending order. The law of ${\bf P}(\alpha,\theta)$ is called the  two-parameter Poisson-Dirichlet distribution and is denoted by
      $PD(\theta,\alpha)$.

      Let $\xi_k,
      k=1,...$ be a sequence of i.i.d. random variables with common
      diffusive distribution $\nu$ on $[0,1]$, i.e., $\nu(x)=0$ for
      every $x$ in $[0,1]$. Set
      \be \label{DIRI1}
      \Xi_{\theta,\alpha, \nu}=\sum_{k=1}^{\infty}P_k(\theta,\alpha)\delta_{\xi_k}.
      \ee
      We call the law of $\Xi_{\theta,\alpha, \nu}$ the two-parameter Dirichlet process, denoted by $Dirichlet(\theta,\alpha,\nu)$.

      If $\alpha=0$ in \rf{GEM1}, we get the well known GEM distribution, the Poisson-Dirichlet distribution, and Dirichlet process
      denoted respectively by
      $GEM(\theta)$, $PD(\theta)$, and $Dirichlet(\theta,\nu)$.

      There is a vast literature on GEM distribution, the Poisson-Dirichlet distribution, and Dirichlet process.
      The areas where they appear include Baysian statistics (\cite{Fer73}), combinatorics (\cite{SheLoy66}), ecology (\cite{FCW43}),
      population genetics (\cite{Ewen04}), and random number theory (\cite{Vershik86}).

      The Poisson-Dirichlet distribution was introduced
      by Kingman \cite{Kingman75} to describe the distribution of gene frequencies in a large neutral
      population at a particular locus. In the population genetics setting it is intimately related to the {\it Ewens sampling
      formula} that describes the distribution of the allelic partition of a sample of size $n$ genes selected from the population.
      The component $P_k(\theta)$ represents the
      proportion of the $k$th most frequent alleles. If $u$ is the
      individual mutation rate and $N_e$ is the effective population size, then the
      parameter $\theta =4N_e u$ is the scaled population mutation rate.

      The GEM distribution can be obtained from the Poisson-Dirichlet distribution through a procedure called {\it size-biased
      sampling}. Here is a brief explanation. Consider a population consisting of individuals of countable
      number of different types labelled $\{1,2,...\}$. Assume that the
      proportion of type $i$ individual in the population is $p_i$. A
      sample is randomly selected from the population and the type of
      the selected individual is denoted by $\sigma(1)$. Next remove all
      individuals of type $\sigma(1)$ from the population and then
      randomly select the second sample. This is repeated to get more
     samples. Denote the type of the $i$th selected sample by
     $\sigma(i)$. Then $(p_{\sigma(1)},p_{\sigma(2)},...)$ is called a
     {\it size-biased permutation} of $(p_1,p_2,...)$. Then sequence $X^{\theta}_k, k=1,2,...$ defined in \rf{GEM1} with $\alpha=0$
     has the same distribution as the
      size-biased permutation of ${\bf P}(\theta)={\bf P}(\theta,0)$. The name GEM distribution is termed by Ewens after
      R.C. Grifffiths, S. Engen and J.W. McCloskey for their contributions to the development of the structure. The Dirichlet process first appeared
      in \cite{Fer73}.

      The literature on the study of Poisson-Dirichlet distribution and Dirichlet process with two parameters is relatively small but is growing rapidly.
      Carlton \cite{Car99} includes detailed calculations of moments and parameter estimations of the two-parameter Poisson-Dirichlet
      distribution. The most comprehensive study of the two-parameter Poisson-Dirichlet
      distribution is carried out in Pitman and Yor
      \cite{PitmanYor97}. In \cite{Der97} and the references therein one can find connections between two-parameter Poisson-Dirichlet distribution
      and models in physics including mean-field spin glasses, random map models,
      fragmentation, and returns of a random walk to origin. The two-parameter Poisson-Dirichlet distribution also found its applications in
      macroeconomics and finance (\cite{Aoki06}).

      The Poisson-Dirichlet distribution and its two-parameter counterpart have many similar structures including
      the urn construction in \cite{Hoppe84} and \cite{FengHoppe98}, GEM representation, sampling formula (\cite{Pitman95}), etc..
      A special feature of the two-parameter Poisson-Dirichlet distribution is included in Pitman
      \cite{Pitman96a} where it is shown that the two-parameter Poisson-Dirichlet distribution is the most general distribution whose
      size-biased permutation has the same distribution as the GEM representation \rf{GEM1}.

      The objective of this paper is to establish large deviation principles (henceforth LDP) for $GEM(\theta,\alpha)$, $PD(\theta,\alpha)$, and
      $Dirichlet(\theta,\alpha,\nu)$ when $\theta$ approaches infinity. Noting that for the one-parameter model, $\theta$ is the scaled
      population mutation rate. For fixed individual mutation rate $u$, large $\theta$ corresponds to large population size.
      In the two parameter setting, we
      no longer have the same explanation. But it can be seen from \rf{GEM1} that for nonzero $\alpha$, large $\theta$ plays a very similar role
      mathematically as
      in the case $\alpha=0$.

      LDP for $Dirichlet(\theta,\nu)$ has been established in \cite{LS87} and \cite{DF01} using different methods. Recently in \cite{DF05}, the LDP is
      established for $PD(\theta)$. From \rf{GEM1}, one can see that for every fixed $k$, the impact of $\alpha$ diminishes as $\theta$ becomes large.
      It is thus reasonable to expect similar LDPs between $GEM(\theta,\nu)$ and $GEM(\theta,\alpha,\nu)$. But in $PD(\theta,\alpha)$ and
      $Dirichlet(\theta,\alpha,\nu)$, every term in
      \rf{GEM1} counts. It is thus reasonable to expect that the LDP for $PD(\theta)$ and $Dirichlet(\theta,\nu)$ are different from the corresponding
      LDPs for $PD(\theta,\alpha)$ and $Dirichlet(\theta,\alpha,\nu)$. But it turns out that the impact of $\alpha$ only appears in the LDP for
      $Dirichlet(\theta,\alpha,\nu)$.

      LDP for $GEM(\theta,\alpha)$ is given in
      Section 2. Using Perman's formula and an inductive structure, we establish the LDP for $PD(\theta,\alpha)$ in Section 3.
      The LDP for $Dirichlet(\theta,\alpha,\nu)$ is established in Section 4 using the subordinator representation in \cite{PitmanYor97} and
      a combination of the methods in \cite{LS87} and \cite{DF01}. Further comments are included in Section 5.

      The reference \cite{DZ98} includes all the terminologies and standard techniques on large deviations
     used in this article. Since the state spaces encountered here are all compact, there is no need to distinguish between
     a rate function and a good rate function.

      \section{LDP for GEM}

      Let $E=[0,1]$, and $E^{\infty}$  be the infinite Cartesian product of $E$. Set
      \[
      {\cal E}=\{(x_1,x_2,...) \in E^{\infty}: \sum_{k=1}^{\infty} x_k\leq 1 \},
      \]
      and consider the map

      \[
      G:E^{\infty}\ra {\cal E}, (u_1,u_2,...)\ra (x_1,x_2,..)
      \]
      with
      \[
      x_1=u_1, x_n =u_n (1-u_1)\cdots(1-u_{n-1}), \ n \geq 2.
      \]

      By a proof similar to that used in Lemma 3.1 in \cite{DF05}, one obtains the following lemma.
      \begin{lemma}\label{gemt1}
      For each $k\geq 1$, the family of the laws of $U_k$ satisfies a LDP on $E$ with speed $\theta$ and rate function
      \be\label{addition4}
      I_1(u)= = \left\{\begin{array}{ll}
      \log\frac{1}{1-u},&  u\in [0,1)\\
      \infty,& u=1
      \end{array}\right.
            \ee
       \end{lemma}

      \begin{theorem}\label{t1}
      The family $ \{GEM(\theta,\alpha): \theta>0, 0<\alpha<1 \} $ satisfies a LDP on ${\cal E}$ with speed $\theta$ and rate function
        \[
      S(x_1,x_2,...)= = \left\{\begin{array}{ll}
      \log\frac{1}{1-\sum_{i=1}^{\infty}x_i},&  \sum_{i=1}^{\infty}x_i <1\\
      \infty,& \mb{else}.
      \end{array}\right.
            \]
       \end{theorem}
      \proof Since $U_1,U_2,..$ are independent, for every fixed $n$ the law of $(U_1,...,U_n)$ satisfies a LDP with
      speed $\theta$ and rate function $\sum_{i=1}^n I_1(u_i)$. For any ${\bf u}, {\bf v}$ in $E^{\infty}$, set
      \[
       |{\bf u}-{\bf v}|=\sum_{i=1}^{\infty}\frac{|u_i-v_i|}{2^i}.
      \]

      Then for any $\delta''>0$ and ${\bf u}$ in $E^{\infty}$, one can choose $n \geq 1$ and small enough
      $0<\delta'<\delta<\delta''$ such that
      \beq
      &&  \{{\bf v}\in E^{\infty}: \max_{1\leq i \leq n}|v_i-u_i|< \delta'\} \subset \{{\bf v}\in E^{\infty}: |{\bf v}-{\bf u}|< \delta\},\\
      &&  \{{\bf v}\in E^{\infty}: |{\bf v}-{\bf u}|\leq \delta\} \subset\{{\bf v}\in E^{\infty}: \max_{1\leq i \leq n}|v_i-u_i|< \delta''\},
      \eeq
      which implies
      \beqn
      &&\lim_{\delta \ra 0}\liminf_{\theta \ra \infty}\frac{1}{\theta}\log P\{|(U_1,U_2,...)-{\bf u}|< \delta\} \geq -\sum_{i=1}^n I_1(u_i),\label{gem-1}\\
      && \lim_{\delta \ra 0}\limsup_{\theta \ra \infty}\frac{1}{\theta}\log P\{|(U_1,U_2,...)-{\bf u}|\leq \delta\} \leq -\sum_{i=1}^n I_1(u_i).\label{gem-2}
      \eeqn

      Since $E^{\infty}$ is compact, by letting $n$ approach infinity in \rf{gem-1} and \rf{gem-2}, it follows that the law of $(U_1,U_2,...)$
      satisfies a LDP with
      speed $\theta$ and rate function $\sum_{i=1}^{\infty} I_1(u_i).$ The effective domain of
      $\sum_{i=1}^{\infty} I_1(u_i)$ is
     \[
      {\cal C}=\{{\bf u}\in E^{\infty}: u_i <1, \sum_{i=1}^{\infty}u_i <\infty\}
     \]
     and
     \be\label{gem-6}
     \sum_{i=1}^{\infty} I_1(u_i)=\log \frac{1}{\Pi_{i=1}^{\infty}(1-u_i)} \ \mb{on}\  {\cal C}.
     \ee

      Since the map $G$ is continuous, it follows from contraction principle and Lemma~\ref{gemt1} that
       the family $ \{GEM(\theta,\alpha): \theta>0, 0< \alpha <1\} $ satisfies a LDP on ${\cal E}$ with speed $\theta$ and rate function
       \be\label{gem-5}
       \inf\{\sum_{i=1}^{\infty} I_1(u_i): u_1=x_1, u_2(1-u_1)=x_2,...\}.
       \ee

For each $1\leq n \leq +\infty$,
\be\label{gem-3}
(1-u_1)\cdots(1-u_{n})=1-\sum_{i=1}^n x_i.
\ee
Hence if $\sum_{i=1}^n x_i=1$ for some finite $n$, then one of $u_1,...,u_n$ is one, and ${\bf u}$ is not in ${\cal C}$.
If $\sum_{i=1}^n x_i <1$ for all finite $n$ and $\sum_{i=1}^{\infty} x_i =1$, then from \rf{gem-3}
\be\label{gem-4}
\lim_{n\ra \infty}\sum_{i=1}^n log(1-u_i)=-\infty,
\ee
which implies that $\sum_{i=1}^{\infty} u_i=\infty$. Thus the inverse of the set $\{(x_1,...)\in {\cal E}:\sum_{i=1}^{\infty}x_i=1\}$
under $G$ is disjoint with ${\cal C}$. Hence the rate function in \rf{gem-5} is the same as $S(x_1,x_2,...)$.

\hfill $\Box$

      \section{LDP for Two-Parameter Poisson-Dirichlet Distribution}

      In this section, we establish the LDP for $PD(\theta,\alpha)$. Recall that ${\bf
       P}(\theta,\alpha)=(P_1(\theta,\alpha),P_2(\theta,\alpha)...,)$ is the random probability
       measure with law $PD(\theta,\alpha)$. When $\alpha=0,$ we will write $${\bf P}(\theta)={\bf
       P}(\theta,0)=(P_1(\theta),P_2(\theta)...,).$$

     \subsection{Perman's Formula}

     For any constant $C >0, \beta >0$, let
      \[
       h(x)=\alpha C x^{-(\alpha +1)}, x >0,
      \]
and
      \[
       c_{\alpha, \beta}= \frac{\Gamma(\beta +1)(C\Gamma(1-\alpha))^{\beta/\alpha}}{\Gamma(\beta/\alpha
       +1)}.
       \]

\noindent
Let $\psi(t)$ be a density function over $(0,\infty)$
      such that for all $\beta >-\alpha$

      \be\label{perman1}
      \int_0^{\infty}t^{-\beta}\psi(t)dt =\frac{1}{c_{\alpha,
      \beta}}.
      \ee
      Set
      \beqn
      && \psi_1(t,p) = h(tp)t \psi(t\bar{p}),
      t>0, 0<p<1, \bar{p}=1-p \label{per1}\\
      && \psi_{n+1}(t,p)=\left\{\begin{array}{ll}
      h(tp) t\int_{p/\bar{p}}^1 \psi_n(t\bar{p}, q)d\,q,& p \leq 1/(n+1)\\
      0,& \mb{else.}
      \end{array}\right.\label{per2}
      \eeqn

      Then the following result is found in \cite{Perman93}(see also \cite{PitmanYor97}).
      \begin{lemma}\label{le-per}
      {\rm (Perman's Formula)} For each $k \geq 1$, let $f(p_1,...,p_k)$ denote the joint density function of
      $(P_1(\alpha,\theta),...,P_k(\alpha,\theta))$. Then
      \be\label{per3}
      f(p_1,...,p_k) =
      c_{\alpha,\theta}\int_0^{\infty}t^{-\theta}g_k(t,p_1,...,p_k)d\,t,
      \ee
      where for $k\geq 2, t>0$, $0<p_k<\cdots<p_1, \sum_{i=1}^k p_i <1$, and $\hat{p}_k=1-p_1-\cdots-p_{k-1}$,
      \be\label{per4}
       g_k(t,p_1,...,p_k)= \frac{t^{k-1}h(tp_1)\cdots
       h(tp_{k-1})}{\hat{p}_k}g_1(t\hat{p}_k, \frac{p_k}{\hat{p}_k})
      \ee
      and
      \be\label{per5}
       g_1(t,p)= \sum_{n=1}^{\infty}(-1)^{n+1}\psi_n(t,p).
      \ee
      \end{lemma}

     \subsection{LDP for $PD(\theta,\alpha)$}

      We will prove the LDP for $PD(\theta,\alpha)$ by first establishing the LDP for $P_1(\theta,\alpha)$ and
      \linebreak $(P_1(\theta,\alpha),...,P_k(\theta,\alpha))$ for any $k\geq 2$.

       \begin{lemma}\label{P1}
      The family of the laws of  $P_1(\theta,\alpha)$ satisfies a LDP on $E$ with speed
      $\theta$ and rate function $I_1(p)$ given in \rf{addition4}.
      %\[
      %I_1(p)= = \left\{\begin{array}{ll}
      %\log\frac{1}{1-p},& p \in [0,1)\\
      %\infty,& \mb{else.}
      %\end{array}\right.
       %     \]
      \end{lemma}
      \proof It follows from the GEM representation that
      \beq
       && E[e^{\lambda \theta U_1}] \leq E[e^{\lambda \theta
       P_1(\theta,\alpha)}]\ \mb{for}\  \lambda \geq 0;\\
       && E[e^{\lambda \theta U_1}] \geq E[e^{\lambda \theta
       P_1(\theta,\alpha)}]\ \mb{for}\  \lambda < 0.
      \eeq
      On the other hand, from the representation in
      Proposition 22 of \cite{PitmanYor97} we obtain that

      \beq
       && E[e^{\lambda \theta P_1(\theta,\alpha)}] \leq E[e^{\lambda \theta
       P_1(\theta)}]\ \mb{for}\  \lambda \geq 0;\\
       && E[e^{\lambda \theta P_1(\theta,\alpha)}] \geq E[e^{\lambda \theta
       P_1(\theta)}]\ \mb{for}\  \lambda < 0.
      \eeq

       Since both the laws of $U_1$ and $P_1(\theta)$ satisfy LDPs with speed $\theta$ and rate
      function $I_1(\cdot)$, we conclude from Lemma 2.4 of
      \cite{DF05} that the law of
      $P_1(\theta,\alpha)$ satisfies a LDP with speed $\theta$ and rate
      function $I_1(\cdot)$.

      \hfill $\Box$

       \begin{lemma}\label{P2}
      For each fixed $k \geq 2$, let
      \[
      \nabla_k =\{{\bf p}=(p_1,...,p_k): 0 \leq p_k\leq\cdots\leq p_1,\sum_{i=1}^k p_i \leq 1
      \},
      \]
      and ${\cal P}_{\theta,k}$ be law of
      $(P_1(\theta,\alpha),...,P_k(\theta,\alpha))$. Then the family $\{{\cal P}_{\theta,k}: \theta >0\}$ satisfies a LDP on $\nabla_k$ with speed
      $\theta$ and rate function
      \[
      I_k(p_1,...,p_k)= = \left\{\begin{array}{ll}
      \log\frac{1}{1-\sum_{i=1}^k p_i},& \sum_{i=1}^k p_i <1\\
      \infty,& \mb{else.}
      \end{array}\right.
            \]
      \end{lemma}
      \proof For $k\geq 2, t>0$, $0<p_k<\cdots<p_1, \sum_{i=1}^k p_i <1$, it follows from \rf{per3}, \rf{per4} and \rf{per5} that
      \beqn
      f(p_1,...,p_k) &=& \frac{c_{\alpha,\theta}(\alpha C)^{k-1}}{\hat{p}_k(p_1\cdots
      p_{k-1})^{\alpha+1}}\int_0^{\infty}t^{-[\theta
      +(k-1)\alpha]}g_1(t \hat{p}_k,
      \frac{p_k}{\hat{p}_k})d\,t\nn\\
      &=& \frac{c_{\alpha,\theta}(\alpha C)^{k-1}\hat{p}_k^{\theta
      +(k-1)\alpha}}{\hat{p}_k^2(p_1\cdots
      p_{k-1})^{\alpha+1}}\int_0^{\infty}s^{-[\theta
      +(k-1)\alpha]}g_1(s,
      \frac{p_k}{\hat{p}_k})d\,s\label{per6}\\
      &=&\frac{c_{\alpha,\theta}(\alpha C)^{k-1}\hat{p}_k^{\theta
      +(k-1)\alpha}}{\hat{p}_k^2(p_1\cdots
      p_{k-1})^{\alpha+1}} \sum_{n=1}^{\infty}(-1)^{n+1}\int_0^{\infty}s^{-[\theta
      +(k-1)\alpha]}\psi_n(s,
      \frac{p_k}{\hat{p}_k})d\,s\nn
      \eeqn

      Set
      \[
      \phi_n(u)=\int_0^{\infty}s^{-[\theta
      +(k-1)\alpha]}\psi_n(s,u)d\,s.
      \]

      Then
      \beqn
      \phi_1(u)&=& \int_0^{\infty}s^{-[\theta +(k-1)\alpha]}s h(s
      u)\psi(s\bar{u})d\,s\nn\\
      &=& (\alpha C)u^{-(\alpha +1)}\int_0^{\infty}s^{-[\theta
      +k\alpha]}\psi(s\bar{u})d\,s\nn\\
      &=& (\alpha C)u^{-(\alpha +1)}\bar{u}^{\theta +k\alpha -1}\int_0^{\infty}s^{-[\theta +k\alpha]}\psi(s)d\,s\label{per7}\\
      &=&(\alpha C)u^{-(\alpha +1)}\bar{u}^{\theta +k\alpha
      -1}\frac{1}{c_{\alpha,\theta +k\alpha}}.\nn
      \eeqn

      For any $n \geq 1$, and $u \leq \frac{1}{n+1}$, it follows from \rf{per2} that

      \beqn
      &&\psi_{n+1}(t,u)\label{per8}\\
      &&\ \ = t\, h(t u)\int_{u/\bar{u}}^1 \psi_{n}(t\bar{u},
      u_1)d\,u_1\nn\\
      &&\ \ = t\, h(t u)\int_{u/\bar{u}}^1 [\chi_{\{u_1\leq 1/n\}}(t\bar{u})h(t\bar{u}u_1)\int_{u_1/\bar{u}_1}^1
      \psi_{n-1}(t\bar{u}\bar{u}_1,
      u_2)d\,u_2]d\,u_1\nn\\
      &&\ \ = t \, h(t u)\int_{u/\bar{u}}^{1/(n-1)}\cdots\int_{u_{n-1}/\bar{u}_{n-1}}^1
      [(t\bar{u}h(t\bar{u}u_1))\cdots (t\bar{u}\bar{u}_1\cdots\bar{u}_{n-2}
      h(t\bar{u}\bar{u}_1\cdots\bar{u}_{n-2}u_{n-1}))]\nn\\
      &&\hspace{1.5cm} \times [\chi_{\{u_1\leq 1/n\}}\cdots \chi_{\{u_{n-1}\leq 1/2\}}]\psi_1(t\bar{u}\bar{u}_1\cdots\bar{u}_{n-2}\bar{u}_{n-1},
      u_n)d\,u_n d\,u_{n-1}\cdots d\,u_1\nn\\
      &&\ \ = t \, h(t u)\int_{u/\bar{u}}^1\cdots\int_{u_{n-1}/\bar{u}_{n-1}}^1
      [(t\bar{u}h(t\bar{u}u_1))\cdots (t\bar{u}\bar{u}_1\cdots\bar{u}_{n-2}
      h(t\bar{u}\bar{u}_1\cdots\bar{u}_{n-2}u_{n-1}))]\nn\\
      && \hspace{1.5cm}\times [\chi_{\{u_1\leq 1/n\}}\cdots \chi_{\{u_{n-1}\leq 1/2\}}][(t\bar{u}\bar{u}_1\cdots\bar{u}_{n-2}\bar{u}_{n-1})
   h(t\bar{u}\bar{u}_1\cdots\bar{u}_{n-2}\bar{u}_{n-1}
      u_n)]\nn\\
      && \hspace{1.5cm}\times\psi(t\bar{u}\bar{u}_1\cdots\bar{u}_{n-2}\bar{u}_{n-1}
      \bar{u}_n)d\,u_n d\,u_{n-1}\cdots d\,u_1\nn\\
      &&\ \ =(\alpha C)^{n+1}t^{-(n+1)\alpha}\int_{u/\bar{u}}^1\int_{u_1/\bar{u}_1}^1\cdots\int_{u_{n-1}/\bar{u}_{n-1}}^1
      [\chi_{\{u_1\leq 1/n\}}\cdots \chi_{\{u_{n-1}\leq 1/2\}}]d\,u_n d\,u_{n-1}\cdots d\,u_1\nn\\
      && \hspace{1.5cm}\times \{(uu_1\cdots u_n)^{-(\alpha
      +1)}\bar{u}^{-n\alpha}\bar{u}_1^{-(n-1)\alpha}\cdots
      \bar{u}_{n-1}^{-\alpha}\psi(t\bar{u}\bar{u}_1\cdots\bar{u}_{n-2}\bar{u}_{n-1}
      \bar{u}_n)\}.\nn
      \eeqn

     Integrating over $t$ we get
     \beqn
     \phi_{n+1}(u)&=&(\alpha C)^{n+1}\int_{u/\bar{u}}^1\int_{u_1/\bar{u}_1}^1\cdots\int_{u_{n-1}/\bar{u}_{n-1}}^1
     [\chi_{\{u_1\leq 1/n\}}\cdots \chi_{\{u_{n-1}\leq 1/2\}}]d\,u_n d\,u_{n-1}\cdots d\,u_1\nn\\
      && \times \{(uu_1\cdots u_n)^{-(\alpha
      +1)}\bar{u}^{-n\alpha}\bar{u}_1^{-(n-1)\alpha}\cdots
      \bar{u}_{n-1}^{-\alpha}\nn\\
     && \times [\int_0^{\infty}t^{-[\theta +(k+n)\alpha]}\psi(t\bar{u}\bar{u}_1\cdots\bar{u}_{n-2}\bar{u}_{n-1}
      \bar{u}_n)d\,t]\}\label{per9}\\
      &=&(\alpha C)^{n+1}u^{-(\alpha+1)}\bar{u}^{(\theta+k\alpha
      -1)}\frac{1}{c_{\alpha, \theta
      +(n+k)\alpha}}A_n(\alpha,\theta)(u)\nn\\
      &=& \phi_1(u)\frac{c_{\alpha,\theta +k\alpha}}{c_{\alpha, \theta
      +(n+k)\alpha}}A_n(\alpha,\theta)(u),\nn
      \eeqn
      where
      \beq
      A_n(\alpha,\theta)(u)&=& \int_{u/\bar{u}}^1\int_{u_1/\bar{u}_1}^1\cdots\int_{u_{n-1}/\bar{u}_{n-1}}^1
      [\chi_{\{u_1\leq 1/n\}}\cdots \chi_{\{u_{n-1}\leq 1/2\}}]d\,u_n d\,u_{n-1}\cdots d\,u_1\\
      && \times \{(u_1\cdots u_n)^{-(\alpha
      +1)}\bar{u}_1^{\theta +(k+1)\alpha-1}\cdots
      \bar{u}_{n}^{\theta +(k+n)\alpha-1}\}.
      \eeq

      Let
      \[
       D = \{(u_1,...,u_n): u_1 \in [u/\bar{u},1/n],...,u_{n-1}\in [u_{n-2}/\bar{u}_{n-1},1/2], u_n \in
       [u_{n-1}/\bar{u}_n,1]\}.
      \]

By definition, $A_n(\alpha,\theta)(u)=0$ for $u=\frac{1}{n+1}$.

For $0< u < 1/(n+1)$, the Lebesgue measure of $D$ is strictly
positive. It follows by direct calculation that

\be\label{per13}
 lim_{\theta \ra \infty}\frac{1}{\theta}\log A_n(\alpha,\theta)(u) \leq \esssup\{\sum_{r=1}^n \log (1-u_r): (u_1,...,u_n)\in D\} <0.
 \ee
On the other hand, by Stirling's formula,
      \be\label{per12}
      lim_{\theta \ra \infty}\frac{1}{\theta}\log\frac{c_{\alpha,\theta}}{c_{\alpha, \theta
      +(n+k)\alpha}}(u)=0.
      \ee
Thus for $n \geq 1$, and $0<u \leq 1/(n+1)$,

      \be\label{per10}
      \lim_{\theta \ra \infty} \frac{c_{\alpha,\theta}}{c_{\alpha, \theta
      +(n+k)\alpha}}A_n(\alpha,\theta) =0.
      \ee

      Let
      \[
       \nabla_k^{\circ}=\{(p_1,...,p_k)\in \nabla_k: p_k >0,
      \sum_{i=1}^k p_i<1\}.
      \]
      Now for each fixed $(p_1,...,p_k)$ in $\nabla_k^{\circ}$, set
      \[
      m= \max\{j\geq 1:   p_k/\hat{p}_k \leq 1/j\},
      \]
      and $A_0(\alpha,\theta)(u)=1$ for any $0<u<1$.

      Then it follows from \rf{per6}, \rf{per9} and \rf{per10}
      that
      \beqn
      f(p_1,...,p_k) &=&\frac{c_{\alpha,\theta}(\alpha C)^{k-1}\hat{p}_k^{\theta
      +(k-1)\alpha}}{\hat{p}_k^2(p_1\cdots
      p_{k-1})^{\alpha+1}} \phi_1(\frac{p_k}{\hat{p}_k})\sum_{n=1}^{m}(-1)^{n+1}\frac{c_{\alpha,\theta +k\alpha}}{c_{\alpha, \theta
      +(n+k-1)\alpha}}A_{n-1}(\alpha,\theta)(\frac{p_k}{\hat{p}_k})\nn\\
      &=&\frac{(\alpha C)^k(\hat{p}_{k+1})^{\theta +k\alpha -1}}{(p_1\cdots p_k)^{1+\alpha}}
      \sum_{n=1}^{m}(-1)^{n+1}\frac{c_{\alpha,\theta}}{c_{\alpha, \theta
      +(n+k-1)\alpha}}A_{n-1}(\alpha,\theta)(\frac{p_k}{\hat{p}_k})\label{per11}\\
      &=&(\hat{p}_{k+1})^{\theta +k\alpha -1}\frac{(\alpha C)^k}{(p_1\cdots p_k)^{1+\alpha}}[1 +o(1)],\nn
      \eeqn
      which implies that
      \be\label{per14}
      \lim_{\theta \ra \infty}\frac{1}{\theta}\log f(p_1,...,p_k)= -\log\frac{1}{1-\sum_{i=1}^k p_i}\ \mb{on}\ \nabla_k^{\circ}.
      \ee

      Introduce a metric $d_k$ on $\nabla_k$ such that for any ${\bf p},{\bf q}$ in $\nabla_k$
      \[
      d_k({\bf p},{\bf q})= \sum_{i=1}^k |p_i-q_i|.
      \]

       For any $\delta >0$, set
       \beq
       V_{\delta}({\bf p})&=&\{{\bf q}\in \nabla_k: d_k({\bf p}, {\bf q})<\delta\}\\
       \bar{V}_{\delta}({\bf p})&=&\{{\bf q}\in \nabla_k: d_k({\bf p}, {\bf q})\leq \delta\}.
       \eeq

       For every ${\bf p}\in \nabla_k^{\circ}$, one can choose $\delta$ small enough such that $V_{\delta}({\bf p})\subset
       \bar{V}_{\delta}({\bf p}) \subset \nabla_k^{\circ}$.
       Let $\mu$ denote the Lebesgue measure on $\nabla_k$. Then by Jensen' inequality and \rf{per14},
       \beqn
       \lim_{\theta\ra \infty}\frac{1}{\theta}\log {\cal P}_{\theta,k}\{V_{\delta}({\bf p})\}
       &=& \lim_{\theta\ra \infty}\frac{1}{\theta}\log\frac{\mu(V_{\delta}({\bf p}))}{\mu(V_{\delta}({\bf p}))}
       \int_{V_{\delta}({\bf p})}f({\bf q})\mu(d\,{\bf q})\label{per15}\\
       &\geq& -\frac{1}{\mu(V_{\delta}({\bf p}))}\int_{V_{\delta}({\bf p})}I_k({\bf q})\mu(d\,{\bf q}).\nn
       \eeqn
       Letting $\delta$ approach zero and using the continuity of $I_k(\cdot)$ at ${\bf p}$, one gets
       \be\label{per16}
       \lim_{\delta \ra 0}\lim_{\theta\ra \infty}\frac{1}{\theta}\log {\cal P}_{\theta,k}\{V_{\delta}({\bf p})\}\geq -I_k({\bf p}).
       \ee

       Since the family $\{{\cal P}_{\theta,k}: \theta >0\}$ is exponentially tight, a partial LDP holds (\cite{Puk91}). Let $J$ be any rate function
       associated with certain subsequence of $\{{\cal P}_{\theta,k}: \theta >0\}$. Then it follows from \rf{per16} that for any ${\bf p}$ in
       $\nabla_k^{\circ}$
       \be\label{per17}
        J({\bf p}) \leq I_k({\bf p}).
       \ee
       Because of the continuity of $I_k$ and the lower semi-continuity of $J$, \rf{per17} holds on $\nabla_k$.

       On the other hand for any ${\bf p}$ in
       $\nabla_k^{\circ}$,
       \beqn
       \lim_{\theta\ra \infty}\frac{1}{\theta}\log {\cal P}_{\theta,k}\{\bar{V}_{\delta}({\bf p})\}
       &=& \lim_{\theta\ra \infty}\frac{1}{\theta}\log \int_{\bar{V}_{\delta}({\bf p})}f({\bf q})\mu(d\,{\bf q})\label{per18}\\
       &\leq& \lim_{\theta\ra \infty}\frac{1}{\theta}\log f({\bf q}_{\delta})\nn\\
       &=& -I_k({\bf q}_{\delta}),\nn
       \eeqn
       where ${\bf q}_{\delta}$ is in $\nabla_k^{\circ}$ such that
       \[
        f({\bf q}_{\delta})=\sup\{f({\bf q}): {\bf q}\in \bar{V}_{\delta}({\bf p}) \}.
       \]
        The existence of such ${\bf q}_{\delta}$ is due to the continuity of $f$ over $\nabla_k^{\circ}$. Letting $\delta$ approach zero, one has

       \be\label{per19}
       \lim_{\delta \ra 0}\lim_{\theta\ra \infty}\frac{1}{\theta}\log {\cal P}_{\theta,k}\{\bar{V}_{\delta}({\bf p})\}\leq -I_k({\bf p}).
       \ee

       Next consider the case that ${\bf p}$ is such that $p_k >0, \sum_{i=1}^k p_i =1$. Then $p_k/\hat{p}_k=1$. For small enough $\delta$, we have
       \[
        q_k/\hat{q}_k > 1/2\ \ \mb{for}\ \ {\bf q} \in \bar{V}_{\delta}({\bf p}).
       \]
       Thus $A_n(\alpha,\theta)(u)=0$ for all $n \geq 1$ on $\bar{V}_{\delta}({\bf p})$ and it follows from \rf{per11} that

       \beqn
       \lim_{\theta\ra \infty}\frac{1}{\theta}\log {\cal P}_{\theta,k}\{\bar{V}_{\delta}({\bf p})\}
       &=& \lim_{\theta\ra \infty}\frac{1}{\theta}\log \int_{\bar{V}_{\delta}({\bf p})}f({\bf q})\mu(d\,{\bf q})\label{per20}\\
       &\leq& \lim_{\theta\ra \infty}\frac{1}{\theta}\log \int_{\bar{V}_{\delta}({\bf p})}(\hat{q}_{k+1})^{\theta +k\alpha -1}\frac{(\alpha C)^k}{(q_1\cdots q_k
      )^{1+\alpha}}\mu(d\,{\bf q}),\nn\\
      &\leq & \log (1-a_{\delta}),\nn
       \eeqn
       where
        $a_{\delta}$ is such that
       \[
        a_{\delta}=\inf\{\sum_{i=1}^k q_i: {\bf q}\in \bar{V}_{\delta}({\bf p})\} <1.
       \]
       Letting $\delta$ go to zero, one gets
       \be\label{per21}
       \lim_{\delta \ra 0}\lim_{\theta\ra \infty}\frac{1}{\theta}\log {\cal P}_{\theta,k}\{\bar{V}_{\delta}({\bf p})\}\leq -I_k({\bf p}).
       \ee

        The only case remains is when there is a $l\leq k$ such that $p_l =0$. The upper bound in this case is obtained by focusing on a lower
        dimensional space of the positive coordinates.

        Thus we have shown that for every ${\bf p}$ in $\nabla_k$
        \be\label{per22}
\lim_{\delta \ra 0}\lim_{\theta\ra \infty}\frac{1}{\theta}\log {\cal P}_{\theta,k}\{V_{\delta}({\bf p})\}=
\lim_{\delta \ra 0}\lim_{\theta\ra \infty}\frac{1}{\theta}\log {\cal P}_{\theta,k}\{\bar{V}_{\delta}({\bf p})\}= -I_k({\bf p}),
        \ee
which combined with the exponential tightness implies the result.

\hfill $\Box$

Let \[
      \bar{\nabla}=\{(p_1,p_2,... ): p_1 \geq p_2\geq \cdots \geq 0,
      \sum_{i=1}^{\infty}p_i\leq 1\}.
      \]
      and for notational simplicity we use ${\cal P}_{\theta}$ to denote the law of ${\bf P}(\alpha,\theta)$ on $\bar{\nabla}$ in the next theorem.
      Then we have
\begin{theorem}\label{P3}
      The family $\{{\cal P}_{\theta}:\theta >0\}$ satisfies a LDP with speed
      $\theta$ and rate function

      \be\label{rate function3}
      I({\bf p})=
      \left\{\begin{array}{ll}
      \log\frac{1}{1-\sum_{i=1}^{\infty} p_i},& (p_1,p_2,...)\in \bar{\nabla},\; \sum_{i=1}^{\infty}p_i < 1 \\
      \infty,& \mb{else.}
      \end{array}\right.
      \ee
      \end{theorem}
      \proof Because $\bar{\nabla}$ is compact, the family $\{{\cal P}_{\theta}:\theta >0\}$ is exponentially tight. It is thus sufficient to verify
      the local LDP (\cite{Puk91}). The topology on
      $\bar{\nabla}$ can be generated by the following metric
      \[
      d({\bf p}, {\bf q})=\sum_{i=1}^{\infty}\frac{|p_i-q_i|}{2^i},
      \]
      where ${\bf p}=(p_1,p_2,...), {\bf q}=(q_1,q_2,...)$. For any fixed $\delta >0$, let $B({\bf
      p},\delta)$ and $\bar{B}({\bf p},\delta)$ denote the respective
      open and closed balls centered at ${\bf p}$ with radius $\delta
      >0$.  Set $n_{\delta}=2+
      [\log_2(1/\delta)]$ where $[x]$ denotes the integer part of $x$.
      Set
      \beq
      &&V_{n_{\delta}}({\bf p};\delta/2)=\{(q_1,q_2,...) \in \bar{\nabla}:
      |q_i-p_i|< \delta/2, i=1,...,n_{\delta}\},\\
      &&V((p_1,...,p_{n_{\delta}});\delta/2)=\{(q_1,...,q_{n_{\delta}}) \in \nabla_{n_{\delta}}:
      |q_i-p_i|< \delta/2, i=1,...,n_{\delta}\}.
      \eeq
      Then we have
      \[
      V_{n_{\delta}}({\bf p};\delta/2) \subset B({\bf p},\delta) .
      \]
      By lemma~\ref{P2} and the fact that
      \[
      {\cal P}_{\theta}\{V_{n_{\delta}}({\bf p};\delta/2)\} =
      {\cal P}_{\theta,n_{\delta}}\{V((p_1,...,p_{n_{\delta}});\delta/2)\},
      \]
      we get that
      \beqn
      \liminf_{\theta \ra \infty}\frac{1}{\theta}\log {\cal P}_{\theta}\{B({\bf
      p},\delta) \}
      &\geq& \liminf_{\theta \ra \infty}\frac{1}{\theta}\log
      {\cal P}_{\theta,n_{\delta}}
      \{ V((p_1,...,p_{n_{\delta}});\delta/2)\} \label{lb1}\\
      &\geq& -I_{n_{\delta}}(p_1,...,p_{n_{\delta}})\geq -I({\bf p})\nn.
      \eeqn
      On the other hand for any fixed $n \geq 1,\delta_1 >0$, let
      \beq
      && U_n({\bf p};\delta_1)=\{(q_1,q_2,...) \in \bar{\nabla}:
      |q_i-p_i|\leq \delta_1, i=1,...,n\},\\
      &&U((p_1,...,p_{n});\delta_1)=\{(q_1,...,q_{n}) \in \nabla_{n}:
      |q_i-p_i|\leq \delta_1, i=1,...,n\}
      \eeq
      Then we have
      \[
      {\cal P}_{\theta}\{U_{n}({\bf p};\delta_1)\}
      ={\cal P}_{\theta,n}\{U((p_1,...,p_{n});\delta_1)\},
      \]
      and, for $\delta$ small enough,
      \[
      \bar{B}({\bf p},\delta)\subset U_n({\bf p};\delta_1),
      \]
      which implies that
      \beqn
       \lim_{\delta \ra 0}\limsup_{\theta \ra \infty}\frac{1}{\theta}\log
       {\cal P}_{\theta}\{\bar{B}({\bf p},\delta)\} &\leq& \limsup_{\theta \ra \infty}
       \frac{1}{\theta}\log {\cal P}_{\theta,n}\{U((p_1,...,p_n),
      \delta_1)\}\label{ub1}\\
      &\leq& -\inf\{I_n(q_1,...,q_n):(q_1,...,q_n)\in U((p_1,...,p_n),
      \delta_1)\} .\nn
      \eeqn

      Letting $\delta_1$ go to zero, and then $n$
      go to infinity, we get

      \be\label{ub2}
      \lim_{\delta \ra
      0}\limsup_{\theta \ra \infty}\frac{1}{\theta}\log
      {\cal P}_{\theta}\{\bar{B}({\bf p},\delta)\} \leq -I({\bf p}),
      \ee
      which
      combined with \rf{lb1} implies the result.

      \hfill $\Box$

      \section{LDP for Two-Parameter Dirichlet Process}

       Let $M_1(E)$ denote the space of all probability measures on $E$ equipped with the weak topology.
       For any diffusive $\nu$ in $M_1(E)$ with support $E$, let $\xi_1,\xi_2,..$
      be independent and
       identically distributed with common distribution $\nu$. Let $\Xi_{\theta,\alpha,\nu}$ be the two-parameter Dirichlet process defined
       in \rf{DIRI1}.
       %\be\label{twodir01}
      % =\sum_{i=1}^{\infty}P_i(\alpha,\theta)\delta_{\xi_i}.
      % \ee

       Let $\{\sigma(t): t \geq 0, \sigma_0=0\}$ be a subordinator with L\'evy measure $x^{-(1+\alpha)}e^{-x}d\,x$, $x>0$, and $\{\tau(t): t \geq 0,
       \tau_0=0\}$ be a gamma subordinator that is independent of $\{\sigma_t: t \geq 0, \sigma_0=0\}$ and has L\'evy measure $x^{-1}e^{-x}d\,x$, $x>0$.

      \begin{lemma}\label{twodir-t0}
      {\rm (Pitman and Yor)} Let
      \be\label{twodir02}
      \gamma(\alpha,\theta)=\frac{\alpha \tau(\frac{\theta}{\alpha})}{\Gamma(1-\alpha)}.
      \ee
      For each $n \geq 1$, and each partition $0<t_1 <\cdots<t_{n}=1$ of $E$, let $A_i=(t_{i-1},t_{i}]$ for $i=2,...,n$, $A_{1}=[0,t_1]$, and $a_j=\nu(A_j)$.
      Set
     \[
      Y_{\alpha,\theta}(t)=\sigma(\gamma(\alpha,\theta)t), t \geq 0.
      \]
      Then
      the distribution of $(\Xi_{\theta,\alpha,\nu}(A_1),...,\Xi_{\theta,\alpha,\nu}(A_{n}))$ is the same as the distribution of
      $$(\frac{Y_{\alpha,\theta}(a_1)}{Y_{\alpha,\theta}(1)},...,\frac{Y_{\alpha,\theta}(\sum_{j=1}^{n}a_j )
      -Y_{\alpha,\theta}(\sum_{j=1}^{n-1}a_j)}{Y_{\alpha,\theta}(1)}).$$
      \end{lemma}
      \proof Proposition 21 in \cite{PitmanYor97} gives the subordinator representation for $PD(\theta,\alpha)$. The lemma follows this
      representation and the construction outlined on page 254
      in \cite{Pitman96}.

   \hfill $\Box$
\vspace{0.5cm}

      \noindent
      Let
      \beq
      && Z_{\alpha,\theta}(t)=\frac{Y_{\alpha,\theta}(t)}{\theta},\\
      && {\bf Z}_{\alpha,\theta}(t_1,...,t_n)=(Z_{\alpha,\theta}(a_1),...,Z_{\alpha,\theta}(\sum_{j=1}^{n}a_j )
      -Z_{\alpha,\theta}(\sum_{j=1}^{n-1}a_j)).
      \eeq

      By direct calculation, one has

      \beqn
      \psi(\la)&=&\log\phi(\la)=\log E[e^{\la \sigma(1)}] =\int_{0}^{\infty}(e^{\la x}-1) x^{-(\alpha +1)}e^{-x}d\,x\label{twodir1}\\
       &=&\left\{\begin{array}{ll}
      \frac{\Gamma(1-\alpha)}{\alpha}[1-(1-\la)^{\alpha}],& \la\leq 1\\
      \infty,& \mb{else.}
      \end{array}\right.\nn
      \eeqn
      and

      \beqn
      L(\la)&=&\lim_{\theta \ra \infty}\frac{1}{\theta}\log E[e^{ \la \tau(\theta)}]\label{twodir2}\\
      &=&\left\{\begin{array}{ll}
      \log(\frac{1}{1-\la}),& \la<1\\
      \infty,& \mb{else.}
      \end{array}\right.\nn
      \eeqn

       For any
       real numbers $\la_1,...\la_n$, let $\vec{\la}=(\la_1,...,\la_n)$. Then by direct calculation
       \beqn
        \frac{1}{\theta}\log E[\exp\{\theta\langle \vec{\la}, {\bf Z}_{\alpha,\theta}(t_1,...,t_n)\rangle\}]&=&
        \frac{1}{\theta}\log E[\Pi_{i=1}^n (E_{\tau(\theta/\alpha)}
        [\exp\{\la_i  \sigma(1)\}]^{\frac{\alpha a_i}{\Gamma(1-\alpha)}\tau(\theta/\alpha)})]\nn\\
       &=& \frac{1}{\theta}\log E[\exp\{(\sum_{i=1}^n \frac{\alpha\nu(A_i)}{\Gamma(1-\alpha)}\psi(\la_i))\tau(\frac{\theta}{\alpha})\}]\label{twodir6}\\
       &\ra&  \Lambda(\la_1,...,\la_n)=\frac{1}{\alpha} L(\frac{\alpha}{\Gamma(1-\alpha)}\sum_{i=1}^n \nu(A_i)\psi(\la_i)).\nn
       \eeqn

       For $(y_1,...,y_n)$ in $R_{+}^n$, set
       \beqn
       J_{t_1,..,t_n}(y_1,...,y_n)&=&\sup_{\la_1,...,\la_n}\{\sum_{i=1}^n \la_i y_i - \Lambda(\la_1,...,\la_n)\}\label{twodir4}\\
       &=& \sup_{\la_1,...,\la_n  \in (-\infty,1]^n}\{\sum_{i=1}^n \la_i y_i + \frac{1}{\alpha}\log [\sum_{i=1}^n \nu(A_i)(1-\la_i)^{\alpha}]\}\nn
       \eeqn

       \begin{theorem}\label{twodir-t1}
       The family of the laws of ${\bf Z}_{\alpha,\theta}(t_1,...,t_n)$
       on space $R_{+}^n$ satisfies a LDP with speed $\theta$ and rate function \rf{twodir4}.
       \end{theorem}
       \proof First note that both function $\psi$ and function $L$ are essentially smooth. Let
       \[
       {\cal D}_{\Lambda}=\{(\la_1,...,\la_n): \Lambda(\la_1,...,\la_n)<\infty\}, \  {\cal D}^{\circ}_{\Lambda}=\mb{interior of}\  {\cal D}_{\Lambda}.
       \]
       It follows from \rf{twodir2} and \rf{twodir6} that
       \[
       {\cal D}_{\Lambda}=\{(\la_1,...,\la_n):\sum_{i=1}^n \nu(A_i)\frac{\alpha}{\Gamma(1-\alpha)}\psi(\la_i))<1\}.
       \]

       %For $\alpha$ between zero and one,
       %\beq
       % \frac{\Gamma(1-\alpha)}{\alpha}&=& \frac{\Gamma(3-\alpha)}{\alpha(1-\alpha)(2-\alpha)}\\
       % &\geq& \frac{1}{\alpha(1-\alpha)(2-\alpha)} \geq 2,
       %\eeq
       %where $\Gamma(3-\alpha)\geq 1$ follows from H\"{older's inequality}.

       \noindent
       The fact that $\nu$ has support $E$ implies that $\nu(A_i)>0$ for $i=1,...,n$,
       and
        \beq
       {\cal D}_{\Lambda}&=&\{(\la_1,...,\la_n):\sum_{i=1}^n \nu(A_i)[1-(1-\la_i)^{\alpha}]<1\}\\
       &=&\{(\la_1,...,\la_n): \la_i \leq 1, i=1,...,n \}\setminus \{(1,...,1)\}, \\
        {\cal D}^{\circ}_{\Lambda}&=&\{(\la_1,...,\la_n): \la_i < 1, i=1,...,n \}.
       \eeq

       \noindent
       Clearly the function $\Lambda$ is
       differentiable on ${\cal D}^{\circ}_{\Lambda}$ and
       $$grad(\Lambda)(\la_1,...,\la_n)
       =\frac{1}{\Gamma(1-\alpha)}L'(\frac{\alpha}{\Gamma(1-\alpha)}
       \sum_{i=1}^n \nu(A_i)\psi(\la_i))(\nu(A_1)\psi'(\la_1),...,\nu(A_n)\psi'(\la_n)).$$

       A sequence $\vec{\la}_m$ approaches the boundary of ${\cal D}^{\circ}_{\Lambda}$ from inside implies that at least one coordinate sequence
       approaches one. Since the interior of $\{\la: \psi(\la)<\infty\}$ is $(-\infty, 1)$ and $\psi$ is essentially smooth, it follows that
       $\Lambda$ is steep and thus
        essentially smooth.  The theorem then follows from
       G\"{a}rtner-Ellis theorem (\cite{DZ98}).

       \hfill $\Box$

      For $(y_1,...,y_n)$ in $R_{+}^n$ and $(x_1,...,x_n)$ in $E^n$, define
      \[
       F(y_1,..,y_n)= \left\{\begin{array}{ll}
      \frac{1}{\sum_{k=1}^n y_k}(y_1,...,y_n),& \sum_{k=1}^n y_k >0\\
      (0,...,0),& (y_1,...,y_n)=(0,...,0)
      \end{array}\right.
      \]
      and
      \be\label{twodir7}
       I_{t_1,..,t_n}(x_1,...,x_n)=\inf\{J_{t_1,..,t_n}(y_1,...,y_n): F(y_1,...,y_n)=(x_1,...,x_n)\}.
       \ee

      Clearly $I_{t_1,...,t_n}(x_1,...,x_n)=+\infty$ if $\sum_{k=1}^n x_k $ is not one. For $(x_1,...,x_n)$ satisfying $\sum_{k=1}^n x_k=1$, we have
      \beqn
       I_{t_1,..,t_n}(x_1,...,x_n)&=&\inf\{J_{t_1,..,t_n}(a x_1,...,a x_n): a=\sum_{k=1}^n y_k >0\}\label{twodir8}\\
       &=&\inf\{\sup_{(\la_1,...,\la_n)\in (-\infty,1]^n}\{a\sum_{i=1}^n \la_i x_i +\frac{1}{\alpha}\log [\sum_{i=1}^n \nu(A_i)(1-\la_i)^{\alpha}\}:
       a >0\}\nn\\
       &=&\inf\{\sup_{(\la_1,...,\la_n)\in (-\infty,1]^n}\{a-\log a \nn\\
       && -\sum_{i=1}^n a(1-\la_i) x_i +\frac{1}{\alpha}\log [\sum_{i=1}^n \nu(A_i)[a(1-\la_i)]^{\alpha}]\}: a>0\}\nn\\
       &=&\inf\{a-\log a : a >0\}
       +\sup_{(\ga_1,...,\ga_n)\in R_{+}^n}\{\frac{1}{\alpha}\log [\sum_{i=1}^n \nu(A_i)\ga_i^{\alpha}]-\sum_{i=1}^n \ga_i x_i\}\nn\\
       &=&\sup_{(\ga_1,...,\ga_n)\in R_{+}^n}\{\frac{1}{\alpha}\log [\sum_{i=1}^n \nu(A_i)\ga_i^{\alpha}]+1-\sum_{i=1}^n \ga_i x_i\}.\nn
      \eeqn

       \begin{theorem}\label{twodir-t2}
       The family of the laws of $(\Xi_{\theta,\alpha,\nu}(A_1),...,\Xi_{\theta,\alpha,\nu}(A_{n}))$
       on space $E^n$ satisfies a LDP with speed $\theta$ and rate function
       \be\label{twodir9}
       I_{t_1,..,t_n}(x_1,...,x_n)=  \left\{\begin{array}{ll}
      \sup_{(\ga_1,...,\ga_n)\in R_{+}^n}\{\frac{1}{\alpha}\log [\sum_{i=1}^n \nu(A_i)\ga_i^{\alpha}]\\
      \hspace{2cm}+1-\sum_{i=1}^n \ga_i x_i\},& \sum_{k=1}^n x_k =1\\
      \infty,& \mb{else}
      \end{array}\right.
       \ee
       \end{theorem}

       \proof Since $J_{t_1,...,t_n}(0,...,0)=\infty$, the function $F$ is thus continuous on the effective domain of
       $J_{t_1,..,t_n}$. The theorem then follows from Lemma~\ref{twodir-t0} and the contraction principle (remark (c) of Theorem 4.2.1 in \cite{DZ98}).

       \hfill $\Box$

       \noindent
       {\bf Remark.} Let $\Delta_n=\{(x_1,...,x_n)\in E^n: \sum_{i=1}^n x_i =1\}.$ Then the result in Theorem~\ref{twodir-t2} holds with $E^n$
       being replaced by $\Delta_n$.

      Let $B_b(E)$ and $C_b(E)$ denote the sets of bounded measurable functions, and bounded continuous functions on $E$ respectively. For each $\mu$
      in $M_1(E)$, set
      \beqn
      I^{\alpha}(\mu)&=& \sup_{f \geq 0, f \in C_b(E)}\{\frac{1}{\alpha}\log(\int (f(x))^{\alpha} \nu(d\,x)) +1-\int f(x)\mu(d\,x)\},
      \label{ratefunction4}\\
      I^{0}(\mu)&=& \sup_{f \geq 0, f \in B_b(E)}\{\int \log f(x) \nu(d\,x) +1-\int f(x)\mu(d\,x)\}, \label{ratefunction5}\\
      H(\nu|\mu)&=&\sup_{g \in C_b(E)}\{\int g(x)\nu(d\,x)- \log \int e^{g(x)}\mu(d\,x)\}\label{ratefunction46}\\
      &=&\sup_{g \in B_b(E)}\{\int g(x)\nu(d\,x)- \log \int e^{g(x)}\mu(d\,x)\},\nn
      \eeqn
      where $H(\nu|\mu)$ is the relative entropy of $\nu$ with respect to $\mu$.

      \begin{lemma}\label{twodir-t3}
       For any $\mu$ in $M_1(E)$,
       \be\label{twodir10}
       I^{\alpha}(\mu)= \sup\{I_{t_1,..,t_n}(\mu(A_1),...,\mu(A_n)): 0<t_1<t_2<\cdots<t_n=1; n=1,2,...\}.
       \ee
       \end{lemma}
       \proof It follows from Tietze's continuous extension theorem and Luzin's Theorem that we can replace $C_b(E)$
       with $B_b(E)$ in the definition of $I^{\alpha}$. This implies that
       \[
        I^{\alpha}(\mu)\geq \sup\{I_{t_1,..,t_n}(\mu(A_1),...,\mu(A_n)): 0<t_1<t_2<\cdots<t_n=1; n=1,2,...\}.
       \]

       On the other hand, for each nonnegative $f$ in $C_b(E)$, let

       \[
       t_i= \frac{i}{n},\ga_i =f(t_i), i=1,...,n.
       \]

       Then
       \beq
       &&\frac{1}{\alpha}\log(\int (f(x))^{\alpha} \nu(d\,x))-\int f(x)\mu(d\,x)
       =\lim_{n \ra \infty}\{\frac{1}{\alpha}\log [\sum_{i=1}^n \nu(A_i)\ga_i^{\alpha}]-\sum_{i=1}^n \ga_i \mu(A_i)\}\\
       && \leq \sup\{I_{t_1,..,t_n}(\mu(A_1),...,\mu(A_n)): 0<t_1<t_2<\cdots<t_n=1; n=1,2,...\},
       \eeq
       which implies
       \[
        I^{\alpha}(\mu)\leq \sup\{I_{t_1,..,t_n}(\mu(A_1),...,\mu(A_n)): 0<t_1<t_2<\cdots<t_n=1; n=1,2,...\}.
       \]

       \hfill $\Box$

       \noindent
       {\bf Remarks.} It follows from the proof of Lemma~\ref{twodir-t3} that the supremum in \rf{twodir10} can be taken
       over all partitions with $t_1,...,t_{n-1}$
       being the continuity points
       of $\mu$. By monotonically approximating nonnegative $f(x)$ with strictly positive functions from above, it follows from
       the monotone convergence theorem that
       the supremum in both \rf{ratefunction4}
       and \rf{ratefunction5} can be taken over strictly positive bounded functions, i.e.,
       \beqn
      I^{\alpha}(\mu)&=& \sup_{f >0, f \in C_b(E)}\{\frac{1}{\alpha}\log(\int (f(x))^{\alpha} \nu(d\,x)) +1-\int f(x)\mu(d\,x)\},
      \label{addratefunction4}\\
      I^{0}(\mu)&=& \sup_{f > 0, f \in B_b(E)}\{\int \log f(x) \nu(d\,x) +1-\int f(x)\mu(d\,x)\}. \label{addratefunction5}
      \eeqn

      \begin{lemma}\label{twodir-t4}
       \be\label{twodir11}
       I^{0}(\mu)= H(\nu|\mu)
       \ee
      \end{lemma}
      \proof If $\nu$ is not absolutely continuous with respect to $\mu$, then $H(\nu|\mu)=+\infty$. Let $A$ be a set such that
      $\mu(A)=0,\nu(A)>0$ and define
      \[
      f_m(x)=\left\{\begin{array}{ll}
      m,& x \in A\\
      1,& \mb{else}
      \end{array}\right.
      \]

      Then
      \[
      I^0(\mu) \geq \nu(A)\log m \ra \infty \ \mb{as}\ m \ra \infty.
      \]

     Next we assume $\nu \ll \mu $ and denote $\frac{d \nu}{d\mu}(x)$ by $\psi(x)$. By definition,
     $$H(\nu|\mu)=\int \psi(x)\log(\psi(x))\mu(d\,x).$$
     Choosing $f_M(x)=\psi(x)\wedge M$ in the definition of $I^0$. Since the function $x\log x$ is bounded below for non-negative $x$,
     applying the monotone convergence theorem on $\{x: \psi(x)\geq e^{-1}\}$, one gets
     \be\label{twodir13}
      I^0(\mu)\geq \lim_{M \ra \infty}\int f_M(x)\log f_M(x)\mu(d\,x)-\log \int f_M(x)\mu(d\,x) = H(\nu|\mu).
     \ee

     On the other hand, it follows by letting $f(x)=e^{g(x)}$ in the definition of $H(\nu|\mu)$ that
     \[
     H(\nu|\mu)=\sup_{f> 0, f \in B_b(E)}\{\int \log f(x)\nu(d\,x)-\log \int f(x)\mu(d\,x)\}.
     \]
    Since
    \[
    \int f(x)\mu(d\,x)-1 \geq \log \int f(x)\mu(d\,x),
    \]
    we get that
    \be \label{twodir14}
     H(\nu|\mu) \geq  I^0(\mu)
    \ee
    which combined with \rf{twodir13} implies \rf{twodir11}.

    \hfill $\Box$

      \begin{lemma}\label{twodir-t5}
      For any $\mu$ in $M_1(E)$, $0\leq \alpha_1 <\alpha_2<1$,
       \be\label{twodir12}
       I^{\alpha_2}(\mu)\geq I^{\alpha_1}(\mu),
       \ee
       and for any $\alpha$ in $(0,1)$, one can find $\mu$ in $M_1(E)$ satisfying $\nu \ll \mu$ such that
       \be\label{twodir15}
       I^{\alpha}(\mu)> I^{0}(\mu),
       \ee
       \end{lemma}
       \proof By H\"{o}lder's inequality, for any $0< \alpha_1 <\alpha_2 <1$,
      \be\label{twodir16}
      \frac{1}{\alpha_1}\log(\int (f(x))^{\alpha_1} \nu(d\,x)) \leq \frac{1}{\alpha_2}\log(\int (f(x))^{\alpha_2} \nu(d\,x)),
      \ee
       and the inequality becomes strict if $f(x)$ is not constant almost surely under $\nu$. Hence
       $I^{\alpha}(\mu)$ is non-decreasing in $\alpha$ over $(0,1)$. It follows from the concavity of $\log x$ that
       \[
        \frac{1}{\alpha}\log(\int (f(x))^{\alpha} \nu(d\,x)) \geq \int \log f(x)\nu(d\,x),
       \]
       which implies that $I^{\alpha}(\mu)\geq I^0(\mu)$ for $\alpha >0$.

       Next choose $\mu$ in $M_1(E)$ such that $\nu \ll \mu$ and $\frac{d \nu}{d\mu}(x)$ is not a constant with $\nu$ probability one,
       then $I^{\alpha}(\mu)> I^{\alpha/2}(\mu)\geq I^0(\mu)$ for $\alpha >0$.

      %As before, let
      % $f_M(x)=\psi(x)\wedge M$. Then
      % \beqn
      % I^{\alpha}(\mu) &\geq& \lim_{M\ra \infty}[ \frac{1}{\alpha}\log(\int \frac{1}{\alpha}(f_M(x))^{\alpha}\log (f_M(x))^{\alpha} \mu(d\,x))
      % +1-\int f_M(x)\mu(d\,x)] \label{twodir18}\\
      % &=& I^0(\mu)+\frac{1}{\alpha}\log(\int (\psi(x))^{\alpha} \nu(d\,x))-\int \log \psi(x)\nu(d\,x) > I^0(\mu).\nn
      % \eeqn

      \hfill $\Box$

      We now ready to prove the main result of this section.

      \begin{theorem}\label{twodir-t5}
       The family of the laws of $\Xi_{\theta,\alpha,\nu}$
       on space $M_1(E)$ satisfies a LDP with speed $\theta$ and rate function
       $I^{\alpha}(\cdot)$.
       \end{theorem}

    \proof Let $\{f_j(x): j=1,2,... \}$ be a countable dense subset of $C_b(E)$ in the supremum norm.
    The set $\{f_j(x): j=1,2,... \}$ is clearly convergence determining
    on $M_1(E)$. Let $|f_j|=\sup_{x \in E}|f_j(x)|$ and
    \[
    {\cal C}=\{g_j(x)=\frac{f_j(x)}{|f_j|\vee 1}: j=1,...\}.
    \]
    Then ${\cal C}$ is also convergence determining.

     For any $\omega, \mu$ in $M_1(E)$, define
    \be\label{twodir20}
     \rho(\omega,\mu)=\sum_{j=1}^{\infty}\frac{1}{2^j}|\langle \omega, g_j\rangle-\langle\mu,g_j\rangle|.
    \ee
   Then $\rho$ is a metric on $M_1(E)$ and generates the weak topology.

    For any $\delta >0, \mu \in M_1(E)$, let
\[
B(\mu,\delta)=\{\omega \in M_1(E): \rho(\omega,\mu)< \delta\}, \ \ \overline{B}(\mu,\delta)=\{\omega \in M_1(E): \rho(\omega,\mu)\leq \delta\}.
\]

    Since $M_1(E)$ is compact, the family of the laws of $\Xi_{\theta,\alpha}$ is exponentially tight. It thus suffices to show that

    \be\label{twodir19}
 \lim_{\delta \ra
0}\liminf_{\theta \ra \infty}\frac{1}{\theta}\log P\{B(\mu,\delta)\} = \lim_{\delta \ra
0}\limsup_{\theta \ra \infty}\frac{1}{\theta} \log P\{\overline{B}(\mu,\delta)\} =- I^{\alpha}(\mu).
\ee

Choose $m$ large enough, one gets
\be\label{addition1}
\{\omega \in M_1(E):|\langle \omega, g_j \rangle -\langle \mu, g_j \rangle|< \delta/2:j=1,\cdots,m \} \subset B(\nu,\delta).
\ee

Choose a partition $t_1, \cdots, t_n $ such that
\[
 \sup\{|g_j(x)-g_j(y)|: x,y \in A_i,i=1,\cdots,n; j=1,\cdots,m \}< \delta/8.
\]

Choosing $0< \delta_1 < \frac{\delta}{4n}$, and define
\[
 V_{t_1, \cdots, t_n}(\mu, \delta_1)= \{(y_1,...,y_n) \in \Delta_n: |y_i-\mu(A_i)|< \delta_1 , i=1,\cdots,n\}.
\]

For any $\omega$ in $M_1(E)$, let
\[
F(\omega)=(\omega(A_1),...,\omega(A_n)).
\]

If $F(\omega) \in V_{t_1, \cdots, t_n}(\mu, \delta_1)$, then for $j=1,...,m$
\beq
|\langle \omega, g_j \rangle -\langle \mu, g_j \rangle|&=&|\sum_{i=1}^n\int_{A_i}g_j(x)(\omega(dx)-\mu(dx))\\
&<& \frac{\delta}{4} + n\delta_1  <  \delta/2,
\eeq
which implies that
\[
\{\omega \in M_1(E): F(\omega)\in V_{t_1, \cdots, t_n}(\mu, \delta_1)\}\subset
\{\omega \in M_1(E):|\langle \omega, g_j \rangle -\langle \mu, g_j \rangle|< \delta/2:j=1,\cdots,m \}.
\]
This combined with \rf{addition1} implies that
\be\label{addition3}
\{\omega \in M_1(E): F(\omega)\in V_{t_1, \cdots, t_n}(\mu, \delta_1)\subset B(\mu,\delta).
\ee

Since $V_{t_1, \cdots, t_n}(\mu, \delta_1)$ is open in $\Delta_{n}$, it follows from Theorem~\ref{twodir-t2} that
\beqn
 &&\lim_{\delta \ra
0}\liminf_{\theta \ra \infty}\frac{1}{\theta}\log P\{B(\mu,\delta)\}\label{twodir22}\\
&&\hspace{1.5cm}\geq
\lim_{\delta \ra
0}\liminf_{\theta \ra \infty}\frac{1}{\theta}\log P\{\omega \in M_1(E): F(\omega)\in V_{t_1, \cdots, t_n}(\mu, \delta_1)\}\nn\\
&&\hspace{1.5cm}= \lim_{\delta \ra
0}\liminf_{\theta \ra \infty}\frac{1}{\theta}\log P\{(\Xi_{\theta,\alpha,\nu}(A_1),...,\Xi_{\theta,\alpha,\nu}(A_{n}))\in V_{t_1, \cdots, t_n}(\mu, \delta_1)\}\nn\\
&&\hspace{1.5cm}\geq -I_{t_1,..,t_n}(\mu(A_1),...,\mu(A_n))\geq -I^{\alpha}(\mu).\nn
\eeqn

Next we will focus on partitions $t_1,...,t_n$ such that $t_1,...,t_{n-1}$ are continuity points of $\mu$. We denote the collection of all such partitions
by ${\cal P}_{\mu}$. This implies that $F(\omega)$ is continuous
 at $\mu$. Hence for any $\delta_2 >0$, one can choose $\delta >0$ small enough such that
\[
\overline{B}(\mu,\delta) \subset F^{-1}\{V_{t_1, \cdots, t_k}(\mu, \delta_2)\}.
\]

Let
\[
 \overline{V}_{t_1, \cdots, t_k}(\nu, \delta_2)=\{(y_1,...,y_n)\in \Delta_n: |y_i-\mu(A_i)|\leq \delta_2 , i=1,\cdots,n-1\}.
\]

Then we have

\beqn
&& \lim_{\delta\ra 0}\limsup_{\theta \ra \infty}\frac{1}{\theta} \log P\{\overline{B}(\mu,\delta)\} \label{w-4}\\
&& \hspace{2cm}\leq  \limsup_{\theta \ra \infty}\frac{1}{\theta} \log P\{(\Xi_{\theta,\alpha,\nu}(A_1),...,\Xi_{\theta,\alpha,\nu}(A_{n}))\in
\overline{V}_{t_1, \cdots, t_n}(\mu, \delta_2)\}.\nn
\eeqn

By
letting $\delta_2$ go to zero and applying Theorem~\ref{twodir-t2} again, one gets

\be\label{w-5}
\lim_{\delta\ra 0}\limsup_{\theta \ra \infty}\frac{1}{\theta} \log P\{\overline{B}(\mu,\delta)\} \leq -
I_{t_1,\cdots,t_n}(\mu(A_1),...,\mu(A_n)).
\ee

Finally, taking supremum over ${\cal P}_{\mu}$ and taking into account the remark after Lemma~\ref{twodir-t3}, one gets
\be\label{w-6}
\lim_{\delta\ra 0}\limsup_{\ga \ra 0}\frac{1}{\theta} \log P\{\overline{B}(\mu,\delta)\} \leq -I^{\alpha}(\mu),
\ee
which, combined with \rf{twodir22}, implies the theorem.

    \hfill $\Box$
\section{Further Comments}

Our results show that the LDPs for $GEM(\theta,\alpha)$ and $PD(\theta,\alpha)$ have the same rate function.
Since $GEM(\theta,\alpha)$ and $PD(\theta,\alpha)$ differs only by the ordering, one would expect to derive the LDP for one from the LDP for the other.
Unfortunately the ordering operation is not continuous and it is not easy to establish an exponential approximation.
The LDPs for $GEM(\theta,\alpha)$ and $PD(\theta,\alpha)$ also have the same rate function as
the LDPs for $GEM(\theta)$ and $PD(\theta)$. Thus $\alpha$ does not play a role in these LDPs.
This is mainly due to the topology used.
It will
be interesting to investigate the possibility of seeing the role of $\alpha$ through establishing the corresponding LDPs on a stronger topology.

The LDPs for $\Xi_{\theta,\alpha,\nu}$ and $\Xi_{\theta,\nu}$ have respective rate functions $I^{\alpha}(\cdot)$ and $I^0(\cdot)$.
Both $\Xi_{\theta,\alpha,\nu}$ and $\Xi_{\theta,\nu}$ converge to $\nu$ for large $\theta$. When $\theta$ becomes large, each
$P_i(\theta,\alpha)$ is more likely to be small. The introduction of positive $\alpha$ plays a similar role.
Thus the mass in $\Xi_{\theta,\alpha,\nu}$ spreads more evenly than the mass in $\Xi_{\theta,\nu}$. Intuitively $\Xi_{\theta,\alpha,\nu}$
is ``closer" to $\nu$ than $\Xi_{\theta,\nu}$. This observation is made rigorous through the fact
that $I^{\alpha}(\cdot)$ can be strictly bigger than $I^0(\cdot)$. The monotonicity of $I^{\alpha}(\cdot)$ in $\alpha$ shows that $\alpha$ can be used
to measure the relative ``closeness" to $\nu$ among all $\Xi_{\theta,\alpha, \nu}$ for large $\theta$.

The process $Y_{\alpha,\theta}(t)$ is a process with exchangeable increments. One could try to establish a general LDP result for processes with
exchangeable increments and derive the result in Section 4 through contraction principle. The proofs here illustrate most of the procedures needed
for pursuing such a general result from which the LDP for $\Xi_{\theta,\alpha,\nu}$ follows.

%\begin{theorem}\label{twodir-t6}
%       \be\label{twodir23}
%       \lim_{\alpha \ra 0}I^{\alpha}(\mu)= I^0(\mu).
%       \ee
%       \end{theorem}
%\proof For any $\alpha$ in $(0,1)$ and any $k\geq 1$, there exists $$

      \end{document}